\documentclass{amsproc}

\usepackage{amssymb}
\usepackage[utf8]{inputenc}
\usepackage[T1]{fontenc}

\usepackage{graphicx}
\usepackage{amssymb}
 \usepackage{amsthm}

\usepackage{epsfig}

\newtheorem{thm}{Theorem}[section]

\theoremstyle{definition}
\newtheorem{defn}[thm]{Definition}

\newtheorem{lem}[thm]{Lemma}
\newenvironment{pf}{\noindent\textbf{Proof.}\quad}{\hfill{$\Box$}}

\theoremstyle{remark}
\newtheorem{rk}[thm]{Remark}

\newcommand{\f}{\rm{Flux}}
\newcommand{\hr}{\textmd{Ham(M,L)}}
\newcommand{\sr}{\textmd{Symp}_0(M,L)}

\newcommand{\usr}{\widetilde{\textmd{Symp}}_0(M,L)}

\begin{document}

\title{On relative Hamiltonian diffeomorphisms}


%    author one information
\author[A.S. Demir]{Ali Sait Demir\\ \textit{Istanbul Technical University, Department of Mathematics, 34469 Istanbul, Turkey}}
\address{Istanbul Technical University, Department of Mathematics, 34469 Istanbul, Turkey}
\curraddr{}
\email{demira@itu.edu.tr}
\thanks{}

%    The 2010 edition of the Mathematics Subject Classification is
%    now available.  If you are citing a classification from the
%    new scheme, use the following input coding instead.
%\subjclass[2010]{Primary }

\date{}

\begin{abstract}

Let $\hr$ denote the group of Hamiltonian diffeomorphisms on a symplectic manifold $M$, leaving a Lagrangian submanifold $L\subset M$ invariant. In this paper, we show that $\hr$ has the fragmentation property, using relative versions of the techniques developed by Thurston and Banyaga.
\end{abstract}

\maketitle
\address
\section{Introduction}

One of the main concerns in the study of automorphism groups of manifolds is whether the group is simple or perfect. The classical technique of Thurston \cite{Thu1} for showing that the group of $C^{\infty}$ diffeomorphisms of a smooth manifold is simple (and hence perfect) requires two main properties of the group: fragmentation and transitivity. Banyaga\cite{Ba1} adaptes these techniques for the group of symplectomorphisms of a symplectic manifold. In this paper we prove that 
$\hr $: the group of relative Hamiltonian diffeomorphisms on a symplectic manifold $M$ satisfies fragmentation property. Roughly the group consists of Hamiltonian diffeomorphisms on $M$ that leave a fixed Lagrangian submanifold $L\subset M$ invariant, (see Ozan \cite{Oz1} for details).
An automorphism group $G$ on a manifold $M$ is called 2-transitive 
 if for all points  in $M$ $x_1,x_2,y_1,y_2$ with $x_1\neq x_2$ and $y_1\neq y_2$, there exists diffeomorphisms $g_1, g_2$ in $G$ with $g_1(x_1)=y_1 ,\ g_2(x_2)=y_2$ and 
 $supp(g_1) \cap supp(g_2)=\emptyset $. However $\hr$ is far from being transitive since it leaves a submanifold invariant. Indeed this is the main 
 reason for the nonsimplicity of $\hr$. The main result of 
 the paper is the Relative Fragmentation Theorem:

\begin{thm} \label{frag}
Let $\mathcal{U}=(U_j)_{j\in I}$ be an open cover of a compact, connected, symplectic manifold $(M,\omega)$ and $h$ be an element of ${\rm{Ham}}(M,L)$ for a Lagrangian submanifold $L$ of $M$. Then $h$ can be written
$$h=h_1 h_2...h_N,$$
where each $h_i\in {\rm{Ham}}_c(M,L), \ i=1,..,N$ is supported in $U_{j(i)}$ for some $j(i) \in I$. Moreover, if $M$ is compact%or $M$ is not compact and $R_{rel}(\Phi)=0$
, we may choose each $h_i$ such that $R_{U_i,U_i\cap L}(h_i)=0$, where we made the identification $U_{j(i)}:=U_i$.
\end{thm}
The relative Calabi homomorphism $R_{U_i,U_i\cap L}$ will be defined in the next section. ${\rm{Ham}}_c(M,L)$ is the subgroup consisting of compactly supported elements.

\noindent For preliminaries and the techniques on the classical diffeomorphism groups  Banyaga's book \cite{Ba1} can be consulted. In this work we followed his exposition of the 
forementioned techniques of Thurston and Banyaga.

\section{Relative Hamiltonian Diffeomorphisms}

Let $(M^{2n}, \omega)$ be a symplectic manifold, i.e. $\omega$ is a closed 2-form such that $\omega^{n}$ is a volume form on $M$. The group of symplectomorphisms is defined as $${\rm{Symp}}(M,\omega)=\{ \phi \in {\rm{Diff}}^{\infty}(M)\ |\ \phi^*\omega= \omega\}.$$
${\rm{Symp}}(M,\omega)$ is by definition equipped with $C^{\infty}$-topology and as first observed by Weinstein in \cite{Wein} it is locally path connected. Let ${\rm{Symp}}_0(M, \omega)$ denote the path component of $id_M \in {\rm{Symp}}(M, \omega)$. For any $\psi \in {\rm{Symp}}_0(M, \omega)$, let $\psi_t \in {\rm{Symp}}(M, \omega)$ for all $t \in [0,1]$, such that $\psi_0=id_M$ and $\psi_1=\psi$. There exists a unique family of vector fields (corresponding to $\psi_t$)

\begin{equation}
X_t: M \longrightarrow TM\ \ \ \textmd{such that} \ \ \ \frac{d}{dt} \psi_t= X_t \circ \psi_t.
\label{eq1}
\end{equation}
The vector fields $X_t$ corresponding to $\psi_t$ (the notation $\dot{\psi_t}$ is also used commonly to denote $X_t$) are called symplectic since they satisfy $\mathcal{L}_{X_t}\omega =0$, where $\mathcal{L}_{X_t}\omega$ denotes the Lie derivative of the form $\omega$ along the vector field $X_t$. By Cartan's formula
$$\mathcal{L}_{X_t}\omega =i_{X_t}(d \omega)+d(i_{X_t}\omega).$$
Hence $X_t$ is a symplectic vector field if and only if $i_{X_t}\omega$ is closed for all $t$. If moreover $i_{X_t}\omega$ is exact, that is to say
$i_{X_t}\omega = dH_t ,\ H_t: M \rightarrow \mathbb{R}$ a family of smooth functions, then $X_t$ are called Hamiltonian vector fields. In this case
the corresponding diffeomorphism $\psi$ is called a Hamiltonian diffeomorphism and $H_1$ is a Hamiltonian for $\psi$. The Hamiltonian diffeomorphisms
form a group as a subgroup in the identity component of the group of symplectomorphisms, ${\rm{Ham}}(M, \omega)\subseteq {\rm{Symp}}_0(M, \omega)$.

Let $(M,\omega)$ be a connected, closed symplectic manifold, $L$ a Lagrangian submanifold of $M$. Denote by $\sr$ the identity component of the group of symplectomorphisms of $M$ that leave $L$ setwise invariant and by $\usr$ its universal cover. Then the restriction of the flux homomorphism to $\usr$ is a well defined homomorphism onto $H^1(M,L)$ ( see \cite{Oz1}), given by
$$\f(\{ \psi_t \})=\int_0^1[i_{X_t}\omega]dt$$
where $\{\psi_t\}\in \widetilde{{\rm{Symp}}}_0(M,L,\omega)$ denotes the homotopy class of smooth paths $\psi_t \in \sr$ with fixed ends $\psi_0=id, \ \psi_1=\psi$ and $X_t$ is the vector field defined by
$\displaystyle{\frac{d}{dt}\psi_t=X_t\circ \psi_t.}$
Note that since $\psi_t$ leaves $L$ invariant, for any $p\in L$, $X_t(p)\in T_pL$.

\noindent \textbf{Notation:} Let $M$ be a manifold, $L\subset M$ a submanifold.
If $f$ is meant to be a map of $M$ that leave $L$ setwise invariant then we write $f:(M,L) \rightarrow (M,L)$.

Let ${\rm{Ham}}(M,L)\subset {\rm{Symp}}_0(M,L)$ be the subgroup consisting of symplectomorphisms $\psi$ such that there is a Hamiltonian isotopy $\psi_t:(M,L)\rightarrow (M,L)$, $t\in [0,1]$ with $\psi_0=id$, $\psi_1=\psi$; i.e. $\psi_t$ is a Hamiltonian isotopy of $M$ such that $\psi_t(L)=L$ for any $t\in [0,1]$. So if $X_t$ is the vector field associated to $\psi_t$ we have $i_{X_t}\omega=dH_t$ for $H_t:M\rightarrow \mathbb{R}$. Since $L$ is Lagrangian $(w|_{L}=0)$, $H_t$ is locally constant on $L$. We have the following characterization which is analogous to its absolute version (see Thm 10.12 in \cite{McD}).
\begin{thm}(\cite{Oz1})
$\psi\in {\rm{Symp}}_0(M,L)$ is a Hamiltonian symplectomorphism if and only if there exists a symplectic isotopy $\psi_t:[0,1]\rightarrow {\rm{Symp}}_0(M,L)$ such that $\psi_0=id$, $\psi_1=\psi$ and $\f(\{\psi_t\})=0$. Moreover, if $\f(\{\psi_t\})=0$ then $\{\psi_t\}$ is isotopic with fixed end points to a Hamiltonian isotopy through points in ${\rm{Symp}}_0(M,L)$.
\end{thm}
\subsection{Relative Calabi Homomorphism}

The relative version of the Calabi homomorphism is defined by the same formula of its absolute version. Let $(M^{2n},\omega)$ be a noncompact symplectic manifold and $L^n\subset M^{2n}$ be a Lagrangian submanifold. If ${\rm{Ham}}_c(M,L)$ is the group of compactly supported Hamiltonian diffeomorphisms of $M$ that leave $L$ invariant, then
$$\rm{R}:\widetilde{{\rm{Ham}}}_c(M,L)\rightarrow \mathbb{R}$$
$$\{\phi_t\}\longmapsto \int_0^1\int_M H_t\omega^ndt \ ,$$
where $H_t$ is given by $i_{X_t}\omega =dH_t$ and $\frac{d}{dt}\phi_t=X_t \circ \phi_t$,
is the relative Calabi homomorphism. That this homomorphism is a well-defined surjective homomorphism can be proved almost the same as the absolute case (see for example \cite{Ba1} p.$103$).

Similarly, the relative Calabi homomorphism can be defined for compact manifolds. Namely, if $\widetilde{{\rm{Ham}}}_{U,U\cap L}(M,\omega)$ denotes the universal cover of Hamiltonian diffeomorphisms supported in $U$ that leave the Lagrangian submanifold $L$ invariant then
$$\rm{R}_{U,U\cap L}:\widetilde{{\rm{Ham}}}_{U,U\cap L}(M,\omega) \rightarrow \mathbb{R}$$
$$\{\phi_t\}\longmapsto \int_0^1\int_M H_t(\omega)^n dt$$
is again a surjective homomorphism.
\begin{rk}
Let $\rm{R}:\tilde{G}\rightarrow \mathbb{R}$ denote any of the above versions of the Calabi homomorphisms in the universal cover setting. We use the same notation for the induced homomorphisms for the underlying groups. Namely, if $\Lambda$ denotes the image of $\pi_1(G)$ under $\rm{R}$, then
$$\rm{R}:G\rightarrow \mathbb{R}/\Lambda$$ is a well-defined homomorphism.
\end{rk}

\subsection{Relative Weinstein Charts}
Let $\psi\in {\rm{Symp}}_0(M,L)$ be sufficiently $C^1$-close to the identity. Similar to the absolute case, there corresponds a closed $1$-form $\sigma=C(\psi)\in \Omega^1(M)$ defined by $\Psi(graph(\psi))=graph(\sigma)$. Here $\Psi:\mathcal{N}(\Delta)\rightarrow \mathcal{N}(M_0)$ is a fixed symplectomorphism between the tubular neighborhoods of the Lagrangian submanifolds diagonal ($\Delta\subset (M\times M,(-\omega)\oplus \omega)$) and the zero section ($M_0\subset (T^{\ast}M,\omega_{can})$) of the cotangent bundle with $\Psi^{\ast}(\omega_{can})=(-\omega)\oplus \omega$. Note that since $\psi\in {\rm{Symp}}_0(M,L)$ the corresponding $1$-form vanish on $TL$, i.e. $\sigma|_{TqL}=0$ for any $q\in L$. Here $\omega_{can}$ denotes the canonical symplectic form on the cotangent bundle of a smooth manifold. See \cite{McD} for the absolute versions and the details.
As a consequence we have the following due to Ozan:

\begin{lem} (\cite{Oz1})
If $\psi\in {\rm{Symp}}_0(M,L,\omega)$ is sufficiently $C^1$-close to the identity and $\sigma=C(\psi_t)\in \Omega'(M)$ then $\psi\in {\rm{Ham}}(M,L)$ iff $[\sigma]\in \Gamma(M,L)$.
\end{lem}

$\Gamma(M,L)$ is the relative flux group defined as the image of the fundamental group of ${\rm{Symp}}_0(M,L,\omega)$ under the flux homomorphism.
$$\Gamma(M,L)=\widetilde{Flux}(\pi_1({\rm{Symp}}_0(M,L,\omega)))\subseteq H^1(M,L,\mathbb{R}).$$
\begin{defn}
The correspondence $$C:{\rm{Symp}}_0(M,L,\omega)\rightarrow Z^1(M,L)$$
$$h\longmapsto C(h)$$
is called a Weinstein chart of a neighborhood of $id_M\in {\rm{Symp}}_0(M,L,\omega)$ into a neighborhood of zero in the set of closed $1$-forms that vanish on $TL$. The form $C(h)$ is called a (relative) Weinstein form.
\end{defn}
With these definitions in mind we have the following. Compare the absolute version in \cite{Ba1}.

\begin{lem} \label{corfrag}
Let $(M,\omega)$ be a symplectic manifold, $L$ a Lagrangian submanifold. For any $h\in {\rm{Ham}}(M,L)$ there exists finitely many hamiltonian diffeomorphisms $h_i\in {\rm{Ham}}(M,L), \ i=1,..,N$, such that each $h_i$ is close to $id_M$ to be in the domain of the Weinstein chart. Moreover $C(h_i)$ is exact for all $i=1,..,N.$
\end{lem}

\begin{pf}
As the above lemma suggests, every smooth path $\psi_t\in {\rm{Ham}}(M,L)$ is generated by Hamiltonian vector fields.
Let $h_t$ be any isotopy in ${\rm{Ham}}(M,L)$ to the identity such that $\frac{d}{dt}h_t=X_t(h_t)$ where $i_{X_t}\omega=df_t$, $h_0=id_M$, $h_1=h$ and
 $f_t:M\rightarrow \mathbb{R}$ are Hamiltonians for all $t \in [0,1]$. Let $N$ be an integer large enough so that
$$\Phi_t^i=\Bigg[h_{\Big(\frac{N-i}{N}\Big)t}\Bigg]^{-1}h_{\Big(\frac{N-i+1}{N}\Big)t}$$
is in the domain of the Weinstein chart. If we let $h_i=\Phi_1^i$ then we have $h=h_Nh_{N-1}...h_1$. As noted by Ozan in \cite{Oz1} the group $\Gamma(M,L)$ is countable. Therefore any continuous mapping of $[0,1]$ into $\Gamma(M,L)$ must be constant.
Hence $t\longmapsto [C(\Phi_t^i)]$ is constant and thus $[C(\Phi_t^i)]=0$.
\end{pf}

\subsection{The Fragmentation Lemma}

By $B^1(M,L)$ denote the set of exact 1-forms that evaluates zero on $TL$ for a Lagrangian submanifold $L\subset M$. To any smooth function $f:M\rightarrow \mathbb{R}$ that is locally constant on $L$, there is a continuous linear map $\sigma_{rel}:B^1(M,L) \rightarrow C^{\infty}_L(M)$ satisfying $\omega=d(\sigma_{rel}(\omega))$ for all $\omega \in B^1(M,L)$. Then there is a bounded linear functional $\tilde{f}:B^1(M,L)\rightarrow B^1(M,L)$, due to Palamadov \cite{Pal}, given by:
\begin{equation}\label{split}
\tilde{f}(\xi)=d(f\sigma_{rel}(\xi))
\end{equation}
We make use of this construction in the proof of the Fragmentation Theorem:
\begin{pf}(of Theorem \ref{frag})
We use the notation of Lemma \ref{corfrag}. By Lemma \ref{corfrag} any $h\in {\rm{Ham}}(M,L)$ can be written as $h=h_1...h_N$ where each $h_i\in {\rm{Ham}}(M,L)$ is close to $id_M$ to be in the domain $V$ of the Weinstein chart
$$C:V\subset {\rm{Symp}}_0(M,L)\rightarrow C(V)\subset Z^1(M,L)$$
and such that $C(h_i)$ is exact for all $i=1,..,N$.

Start with an open cover $\mathcal{U}=(U_i)_{i\in \mathbb{N}}$ of $M$ and a partition of unity $\{ \lambda_i \}$ subordinate to it. Let $K$ be a compact subset of $M$ containing the support of $h$. Let $\mathcal{U}_k=\{U_0,...,U_N \}$ be a finite subcover for $K$ such that $U_i\cap U_{i+1}\neq \emptyset$. Then consider the functions
$$\mu_0=0 \quad , \quad \mu_j=\sum_{i\leq j}\lambda_i $$
for $j=1,2,...,N$. Note that for any $x\in K$, $\mu_N(x)=1$ and $\mu_i(x)=\mu_{i-1}(x)$ for $x \notin U_i$.

Let $\tilde{\mu_i}$ be defined as in the Equation (\ref{split}). Since this operator is bounded, there is an open neighborhood
$ V_0 \subset V$ of $id\in {\rm{Symp}}_c(M,L)$ with
$$\tilde{\mu_i}(C(h)) \in C(V) \ \mbox{for all } i=1,...,N \ \mbox{and } h \in V_0$$
Assume that $h \in V_0$ and define
$$\psi_i=C^{-1}(\tilde{\mu}_i(C(h))) \in {\rm{Ham}}(M,L).$$
Note that $\psi_{i-1}(x)=\psi_i(x)$ for $x \notin U_i$ since $\mu_{i-1}(x)=\mu_i(x)$ in that case. Therefore $(\psi_{i-1}^{-1} \psi_i)(x)=x$ if
$x \notin U_i$. Hence, $h_i=(\psi_{i-1})^{-1}( \psi_i)$ is supported in $U_i$. On $K$ we have $\mu_N=1,\ \mu_0=0, \ \psi_N=h, $ and $\psi_0=id$.
Therefore
$$h=\psi_N=(\psi_{0}^{-1} \psi_1)(\psi_{1}^{-1} \psi_2)...(\psi_{N-1}^{-1} \psi_N)=h_1h_2...h_N.$$

For the second statement define the isotopies $h_t^i=\psi_{i-1}(t) \psi_i(t)$, where $\psi_i(t)=C^{-1}(t \tilde{\mu_i(C(h))})$. A classical result due to Calabi states that the Lie algebra of locally supported Hamiltonian diffeomorphisms is perfect \cite{Ca1}. Since for each
$t, \ \dot{h_t^i}$ is a Hamiltonian vector field parallel to $L$, we can write $\dot{h_t^i}$ as a sum of commutators. In other words we have
$$\dot{h_t^i}=\sum_{j}[X^{ji}_{t},Y^{ji}_{t}],$$
where $X^{ji}_{t}$ and $Y^{ji}_{t} $are again Hamiltonian vector fields (not necessarily parallel to $L$). By the cut-off lemma below $X^{ji}_{t}$ and $Y^{ji}_{t}$ can
be chosen to vanish outside of an open set whose closure contain $U_i$. If $u_t^i$ is the unique function supported in $U_i$ with
$i_{\dot{h_t^i}} \omega = du_t^i$, then $du_t^i=\sum_{j}{\omega(X^{ji}_{t},Y^{ji}_{t})}$ since both functions above have the same differential and both
have compact supports. Therefore
$$\int_{U_i}u_t^i\omega^n=\int_{M}u_t^i\omega^n=\sum_{j} \int_{}\omega(X^{ji}_{t},Y^{ji}_{t})\omega^n=0$$
implying that
$$\int_0^1 \int_{U_i}u_t^i\omega^n=R_{U_i,U_K\cap L}(h_i)=0.$$
\end{pf}

The cut-off lemma we used in the proof of the fragmentation lemma is as follows.

\begin{lem} \label{cutoff}
Let $\varphi_t \in {\rm{Ham}}(M,L)$ be an isotopy of a smooth symplectic manifold $(M,\omega)$ leaving a Lagrangian submanifold $L$ invariant. Let $F \subset M$
be a closed subset and $U,V \subset M$ open subsets such that $U\subset \overline{U} \subset V$ with $\cup_{t \in [0,1]}\varphi_t(F) \subset U$.
Then there is an isotopy $\overline{\varphi_t} \in {\rm{Symp}}(M,L)$ supported in $V$ and coincides with $\varphi_t$ on $U$.
\end{lem}

\begin{pf}
 We choose a smooth function $\lambda_t(x)=\lambda(x,t)$ which equals to $1$ on $U\times[0,1]$, $0$ outside of $V\times [0,1]$. Let $f_t$ denote the family of Hamiltonians corresponding to $\varphi_t$, i.e. $i_{\dot{\varphi_t}}\omega=df_t$. Define
 $\overline{X}(x,t)=X_{(\lambda_t \cdot f_t)}+\partial/ \partial t$ on $M \times [0,1]$, where $X_{(\lambda_t \cdot f_t)}$ is the Hamiltonian
vector field given by $i_{X_{(\lambda_t \cdot f_t)}}\omega=d(\lambda_t \cdot f_t)$.The desired isotopy is obtained by integrating the vector field $\overline{X}(x,t)$.
\end{pf}

\section{Final Remarks}

\begin{enumerate}
\item ${\rm{Ham}}(M,L)$ is not simple because of the following. Consider the sequence of groups and homomorphisms:
$$0\longrightarrow \rm{Ker}\varphi \longrightarrow {\rm{Ham}}(M,L) \stackrel{\varphi}{\longrightarrow} \rm{Diff}^{\infty}(L) \longrightarrow 0,$$
where $\varphi$ is just restriction to $L$. Therefore $\rm{Ker}\varphi$ consists of Hamiltonian diffeomorphisms of $M$ that are identity when restricted to $L$. Clearly, $\rm{Ker}\varphi$ is a closed subgroup.

\item The remaining parts of the Thurston and Banyaga's proofs fail to work in the relative case. In the absolute case the proof is completed by finding one smooth manifold for which the perfectness is easily shown. In both smooth and symplectic categories this is done through the torus due to a theorem by Herman \cite{Her1}. Unfortunately the underlying KAM theory fails to apply in the relative case and thus whether $\hr$ is perfect remains open.
\end{enumerate}

\bibliographystyle{model1a-num-names}

%\bibliographystyle{amsplain}

\end{document}